\documentclass{article}

\newtheorem{fed}{\textbf{Definition}}[section]

\newtheorem{lemma}[fed]{\textbf{Lemma}}

\newtheorem{prop}[fed]{\textbf{Proposition}}

\usepackage{amssymb,bbm,graphicx,epsfig,psfrag,epic,eepic,latexsym}
\usepackage{amsmath}
\usepackage{mathrsfs}
\usepackage[dvips]{color}

\begin{document}
\title{Nullity bounds for certain Hamiltonian delay equations}
\author{Urs Frauenfelder}
\maketitle

\begin{abstract}
In this paper we introduce a class of Hamilton delay equations which arise as critical points of
an action functional motivated by orbit interactions. We show that the kernel of the Hessian at each
critical point of the action functional satisfies a uniform bound on its dimension. 
\end{abstract}
\section{Introduction}

Given a symplectic manifold $(M,\omega)$ and a smooth function $H \colon M \to \mathbb{R}$ the Hamiltonian
vector field $X_H$ of $H$ is implicitly defined by
$$dH=\omega(\cdot, X_H).$$
It is an old problem to study periodic solution of the Hamiltonian vector field, i.e., $u \in C^\infty(S^1,M)$ where
$S^1=\mathbb{R}/\mathbb{Z}$ is the circle, satisfying
$$\partial_t u(t)=X_H(u(t)), \quad t \in S^1.$$
Periodic solution can as well be interpreted as fixed points. Namely consider the flow of the Hamiltonian vector field
$\phi^t_H$ satisfying
$$\phi^0_H=\mathrm{id}_M, \qquad \frac{d}{dt}\phi^t_H=X_H \circ \phi^t_H.$$
Then $u \in C^\infty(S^1,M)$ is a periodic orbit if and only if
$$\phi^1_H(u(0))=u(0),$$
i.e., $u(0)$ is a fixed point of $\phi^1_H$. 

Alternatively periodic orbits arise as well variationally as critical points of the action functional of classical mechanics. To simplify this discussion we assume that the symplectic manifold is exact, i.e., $\omega=d\lambda$. Then the action functional of classical mechanics is defined by
$$\mathcal{A}_H \colon C^\infty(S^1,M) \to \mathbb{R}, \quad u \mapsto \int u^* \lambda-\int_0^1 H(u(t))dt$$
and its critical points are periodic orbits. 

Suppose that $u$ is a critical point of $\mathcal{A}_H$, i.e., a periodic orbit. Assume further that
$\xi \in \ker \mathcal{H}_{\mathcal{A}_H}(u)$ is an element of the kernel of the Hessian of the functional 
$\mathcal{A}_H$ at its critical point $u$. The tangent space of
the free loop space at $u$ consists of smooth vector fields along $u$. In particular, $\xi$ is a smooth vector field
along $u$. In view of the interpretation of critical points as fixed points, we obtain that
$$\xi(t)= d\phi^t_H(u(0)) \xi(0) \in T_{u(t)}M, \qquad t \in S^1.$$
In particular, $\xi(t)$ is completely determined by $\xi(0) \in T_{u(0)} M$. This means that the \emph{nullity} at $u$, i.e., the dimension of the kernel of the Hessian at $u$, is uniformly bounded by
$$\dim \ker \mathcal{H}_{\mathcal{A}_H} \leq \dim T_{u(0)}M=\dim M.$$

Recently the author with his collaborators started studying Hamiltonian delay equations, see \cite{albers-frauenfelder-schlenk1, albers-frauenfelder-schlenk2, frauenfelder-weber}. While for ``reasonable" Hamiltonian delay equations 
the Hessian can be interpreted as a Fredholm operator of index zero and in particular its kernel is therefore always finite dimensional, a uniform bound on its dimension is quite unlikely in general, since due to the nonlocal nature of Hamiltonian delay equations the interpretation of critical points as fixed points of a flow does not longer hold true. 
The purpose of this note is to introduce a class of Hamiltonian delay equations for which such a uniform bound on the nullity still holds although the problems are in general not local. The author discovered this kind of equations in connection with his study on Helium \cite{frauenfelder}. The idea is to study particles which interact with each other
just through their orbits. This is kind of a semiclassical analogon of the Hartree-Fock method in quantum mechanics, see for instance \cite{bethe-salpeter}. Another instance where this kind of Hamiltonian delay equations shows up is related to the recent interesting paper of Barutello, Ortega, and Verzini \cite{barutello-ortega-verzini}. There the authors discovered a nonlocal regularization of the collisions in Kepler problem working for all energies simultanuously. Most geometric regularizations like the one by Moser \cite{moser} or Levi-Civita \cite{levi-civita} depend on the energy. An exception is the one by Ligon and Schaaf \cite{ligon-schaaf}, which however has some issues with smoothness under perturbations. 

This paper is organised as follows. In Section~\ref{class} we introduce the class of Hamiltonian delay equations and
in Section~\ref{examples} we discuss various examples where they show up. In Section~\ref{results} we state the main results showing how the nullity for these Hamiltonian delay equations is uniformly bounded. In general the bound is
higher than just the dimension of the symplectic manifold, however there are cases where it can be shown that the 
dimension is sufficient as upper bound. In Section~\ref{selfadjoint} and Section~\ref{commcase} we discuss a class of
selfadjoint operators for which we can proof a uniform bound on the dimension of its kernel. In Section~\ref{prA} and
Section~\ref{prB} we proof the nullity bounds for our Hamiltonian delay equations by showing that the Hessian gives rise
to selfadjoint operators studied in the previous two Sections. Finally in Appendix~\ref{symmetries} 
we discuss the symmetries of our Hamiltonian delay equations and show that a solution gives rise to iterated solutions.

For classical periodic orbits the Maslov index behaves as a quasimorphisms under iterations \cite{salamon-zehnder}. 
It would be interesting to now if this property continues to hold as well for some Hamiltonian delay equations. This property of the Maslov index for example plays an important role in the proof of the Conley conjecture \cite{ginzburg}.
The Maslov index is as well an important ingredient in the EBK method in semiclassics \cite{gutzwiller}. The author hopes that this paper triggers some research in these directions for Hamiltonian delay equations. 
\\ \\
\emph{Acknowledgements: } The author acknowledges partial support by DFG grant FR 2637/2-2.

\section{A class of Hamiltonian delay equations}\label{class}

Suppose that $(M,\omega)$ is a symplectic manifold which we assume to be symplectically aspherical, i.e., the symplectic
form integrated over any sphere in $M$ vanishes. This happens for example if the symplectic manifold is exact, due to 
Stokes' Theorem. We consider pairs $\mathcal{F}=(f,H)$ of smooth functions
$$H \colon M \to V, \qquad f \colon W \to \mathbb{R}$$
where $V$ is a finite dimensional vector space and $W \subset V$ an open subset. There is no compatibility required like linearity between the vector space
structure of $V$ and the two functions. However, the vector space structure allows us to canonically extend the
function $H$ to the free loop space of $M$. Namely abbreviating by $S^1=\mathbb{R}/\mathbb{Z}$ the circle
the map $H$ induces a map on the free loop space 
$$\overline{H} \colon C^\infty(S^1,M) \to V, \quad u \mapsto \int_0^1 H(u(t))dt.$$
We abbreviate by
$$\mathcal{L} \subset C^\infty(S^1,M)$$
the component of contractible loops in the free loop space of $M$ and set
$$\mathcal{L}_\mathcal{F}=\big(\overline{H}\big)^{-1}(W) \cap \mathcal{L}$$
the open subset of contractible loops whose mean value under $H$ lies in $W$. Note that we do not require
that $H(u(t)) \in W$ for every $t \in S^1$, if $u \in \mathcal{L}_\mathcal{F}$.
If $u \in \mathcal{L}_\mathcal{F}$ it admits, since it is contractible, a filling disk, i.e.,
a smooth map
$$\bar{u}\colon D \to M,$$
where $D=\{z \in \mathbb{C}: |z| \leq 1\}$ is the closed unit disk, satisfying
$$\bar{u}(e^{2\pi it})=u(t), \quad t \in S^1.$$
The action functional for the pair $\mathcal{F}=(f,H)$
$$\mathcal{A}_{\mathcal{F}} \colon \mathcal{L}_\mathcal{F} \to \mathbb{R}$$
is defined for $u \in \mathcal{L}_\mathcal{F}$ by
$$\mathcal{A}_{\mathcal{F}}(u)=\int_D \bar{u}^* \omega-f(\overline{H}(u)).$$
Suppose that $\lambda \in V^*$, i.e., $\lambda$ lies in the dual vector space of $V$. Then the composition
$$\lambda \circ H \colon M \to \mathbb{R}$$
is a smooth function on $M$. The Hamiltonian vector field of $\lambda \circ H$ is implicitly defined by
$$d(\lambda \circ H)=\omega(\cdot, X_{\lambda \circ H})$$
and we obtain a linear map
$$X_H \colon V^* \to \Gamma(TM), \quad \lambda \mapsto X_{\lambda \circ H}.$$
The following Lemma is straightforward.
\begin{lemma}\label{crit}
Critical points of the action functional $\mathcal{A}_{\mathcal{F}}$ are contractible solutions $u \in C^\infty(S^1,M)$ of the problem
$$\partial_t u(t)=X_H\Big(df(\overline{H}(u))\Big)\big(u(t)\big), \quad t \in S^1$$
whose mean value under $H$ lies in $W$.
\end{lemma}

\section{Examples}\label{examples}

We consider several examples of pairs $\mathcal{F}=(f,H)$ as discussed in the previous Section.
\\ \\
\textbf{Example\,1:} Assume that $f \colon V \to \mathbb{R}$ is linear. Then for $u \in \mathcal{L}$
$$\mathcal{A}_{\mathcal{F}}(u)=\int_D \bar{u}^* \omega-f(\overline{H}(u))=\int_D \bar{u}^* \omega-\overline{fH}(u)$$
which is nothing else than the usual action functional of classical mechanics of the Hamiltonian $fH \colon M \to \mathbb{R}$.
In particular, its critical points are contractible one-periodic orbits of the Hamiltonian vector field of $fH$.
\\ \\
\textbf{Example\,2: } Assume that $V=\mathbb{R}$, so that $H \colon M \to \mathbb{R}$ is just a smooth function on $M$.
Suppose further, that 
$$f \colon \mathbb{R} \to \mathbb{R}, \quad x \mapsto \frac{1}{2}x^2.$$
Then the action functional is given for $u \in \mathcal{L}$ by
$$\mathcal{A}_\mathcal{F}(u)=\int_D \bar{u}^*-\frac{1}{2}\overline{H}^2(u).$$
Its critical points are contractible loops $u \colon S^1 \to M$, solving the problem
$$\partial_t u(t)=\overline{H}(u)X_H\big(u(t)\big),\quad t \in S^1.$$
By preservation of energy we see that $H$ is constant along $u$, so that we have
$$H(u(t))=\overline{H}(u), \quad t \in S^1.$$
If we reparametrise $u$ by $\overline{H}(u)$, we can interpret it as a periodic orbit of $X_H$ of period its energy
$\overline{H}(u)$.
\\ \\
\textbf{Example\,3: } Consider $M=T^*(0,\infty) \times T^*(0,\infty)=(0,\infty)\times \mathbb{R} \times (0,\infty) \times
\mathbb{R}$. 
For loops $(q_1,p_1,q_2,p_2) \in C^\infty(S^1,M)$ satisfying $\overline{q}_2>\overline{q}_1$ consider the functional
$$\mathcal{A}(q_1,p_1,q_2,p_2)=\int p_1 dq_1+\int p_2 dq_2-\int_0^1 \bigg(\frac{p_1^2}{2}+\frac{p_2^2}{2}-\frac{\mu}{q_1}
-\frac{\mu}{q_2}\bigg)dt-\frac{1}{\overline{q}_2-\overline{q}_1},$$
where $\mu>1$. This functional describes two electrons attracted by a nucleus of charge $\mu$ and interacting with
each other by their mean position, see \cite{frauenfelder}. If we choose $V=\mathbb{R}^3$,
$$H \colon M \to \mathbb{R}^3: (q_1,p_1,q_2,p_2) \mapsto \bigg(\frac{p_1^2}{2}+\frac{p_2^2}{2}-\frac{\mu}{q_1}
-\frac{\mu}{q_2}, q_1,q_2\bigg),$$
$$W=\big\{(x_0,x_1,x_2) \in \mathbb{R}^3: x_2>x_1\big\},$$
and
$$f \colon W \to \mathbb{R}, \quad (x_0,x_1,x_2)=x_0+\frac{1}{x_2-x_1},$$
then we have
$$\mathcal{A}=\mathcal{A}_\mathcal{F}.$$
\textbf{Example\,4: } On $T^*\mathbb{C}=\mathbb{C}\times \mathbb{C}$ define
$$H \colon T^*\mathbb{C} \to \mathbb{R}^2, \quad (z,w) \mapsto (|z|^2,|w|^2).$$
For the open subset
$$W=(\mathbb{R}\setminus \{0\}) \times \mathbb{R} \subset \mathbb{R}^2$$
we put
$$f \colon W \to \mathbb{R}, \quad (x_1,x_2) \mapsto \frac{x_2-8}{8x_1}.$$
For the tuple $\mathcal{F}=(f,H)$ the functional $\mathcal{A}_\mathcal{F} \colon \mathcal{L}_\mathcal{F} \to \mathbb{R}$
is then given by
$$\mathcal{A}_\mathcal{F}(z,w)=\int(w_1dz_1+w_2dz_2)-\frac{\int_0^1 |w(t)|^2dt-8}{8\int_0^1 |z(t)|^2dt}.$$
Its critical points are solutions of the problem
$$\left\{\begin{array}{c}
\partial_t z(t)=\frac{w(t)}{4\int_0^1|z(s)|^2ds}\\
\partial_tw(t)=\frac{\big(\int_0^1 |w(s)|^2ds-8\big)z(t)}{4\big(\int_0^1 |z(s)|^2ds\big)^2}
\end{array}\right.$$
for $t \in S^1$. If one replace $w$ with the first derivative of $z$ using the first equation the second equation
gives rise to a problem involving only $z$ but now as well its second derivative
\begin{eqnarray}\label{bov}
\partial^2_t z(t)&=&\frac{\partial_t w(t)}{4\int_0^1|z(s)|^2ds}\\ \nonumber
&=&\frac{\big(\int_0^1 |w(s)|^2ds-8\big)z(t)}{16\big(\int_0^1 |z(s)|^2ds\big)^3}\\ \nonumber
&=&\Bigg(\frac{\int_0^1 |\partial_t z(s)|^2ds}{\int_0^1 |z(s)|^2ds}-\frac{1}{2\big(\int_0^1 |z(s)|^2ds\big)^3}\Bigg)z(t).
\end{eqnarray}
This problem appeared recently in the work by Barutello, Ortega and Verzini \cite[Lemma\,3.13]{barutello-ortega-verzini}.
The interesting aspect of this problem is that if one defines for a given solution $z$ of (\ref{bov})
$$t_z \colon [0,1] \to [0,1], \quad \tau \mapsto \frac{\int_0^\tau |z(s)|^2ds}{\int_0^1 |z(s)|^2ds}$$
and denotes by
$$\tau_z \colon [0,1] \to [0,1]$$
the inverse of $t_z$, then the map
$$x \colon S^1 \to \mathbb{C},\quad t \mapsto z^2(\tau_z(t))$$
is a one-periodic solution of the planar Kepler problem having possibly collisions with the mass at the origin.
\\ \\
We say that a tuple $\mathcal{F}=(f,H)$ is \emph{commuting} if for any $\lambda_1,\lambda_2 \in V^*$ the Hamiltonians
$\lambda_1 \circ H$ and $\lambda_2 \circ H$ Poisson commute in the sense that
\begin{equation}\label{comm}
\{\lambda_1 \circ H, \lambda_2 \circ H\}:=\omega(X_{w_1 \circ H}, X_{w_2 \circ H})=0.
\end{equation}
Instances of commuting pairs are obtained by considering a collection of particles each one modelled by a symplectic manifold whose orbits interact with each other. One might think that each particle has its individual time and therefore the interaction cannot be local but involves its whole orbit. The following Example shows the details of this construction.
\\ \\
\textbf{Example\,5: } Suppose that $m \in \mathbb{N}$ and  for $1 \leq j \leq m$ there is given a symplectically aspherical symplectic manifold
$(M_j,\omega_j)$ together with a smooth function $H_j \colon M_j \to \mathbb{R}$. Then the product symplectic manifold
$$(M,\omega)=(M_1\times \ldots \times M_m, \omega_1 \oplus \ldots \oplus \omega_m)$$
is symplectically aspherical as well and we define
$$H \colon M \to \mathbb{R}^m, \quad (x_1,\ldots,x_m) \to (H_1(x_1),\ldots,H_m(x_m)).$$
Then for any smooth $f \colon W \to \mathbb{R}$, where $W$ is an open subset of $\mathbb{R}^m$ the pair
$(f,H)$ is commuting. 

\section{Statement of the main results}\label{results}

Our first main result is the following uniform bound on the dimension of the kernel of the Hessian $\mathcal{H}_{\mathcal{A}_\mathcal{F}}$ at an arbitrary critical point $u$
of the action funtional $\mathcal{A}_\mathcal{F}$.
\\ \\
\textbf{Theorem\,A: } \emph{Suppose that $\mathcal{F}=(f,H)$ is a pair consisting of smooth functions $H \colon M \to V$
and $f\colon W \to \mathbb{R}$, where $V$ is a real finite dimensional vector space and $W \subset V$ is an open subset.  
Assume further that $u$ is a critical point of $\mathcal{A}_\mathcal{F}$. Then
$$\mathrm{dim}  \ker\big(\mathcal{H}_{\mathcal{A}_\mathcal{F}}(u)\big) \leq \mathrm{dim}(M)+\mathrm{dim}(V).$$}
\\ \\
It the pair is commuting as explained in (\ref{comm}) we obtain a stronger upper bound.
\\ \\
\textbf{Theorem\,B: } \emph{Under the assumptions of Theorem\,A assume in addition that $\mathcal{F}$ is commuting, then
$$\mathrm{dim}  \ker\big(\mathcal{H}_{\mathcal{A}_\mathcal{F}}(u)\big) \leq \mathrm{dim}(M).$$}

\section{A class of selfadjoint operators}\label{selfadjoint}

In this section we introduce a class of selfadjoint operators for which we prove a uniform bound on the dimension of
their kernel. This result is used in the proof of Theorem\,A.
\\ \\
We consider $\mathbb{C}^n$ endowed with its canonical symplectic structure. Namely if we think of 
$\mathbb{C}^n$ as $\mathbb{R}^{2n}$ and let $J \colon \mathbb{R}^{2n} \to \mathbb{R}^{2n}$ be the linear map
obtained from multiplication by $i$, we have
$$\omega(\xi, \eta)=\xi^T J \eta, \qquad \xi, \eta \in \mathbb{R}^{2n}.$$
Note that
$$J=-J^T=-J^{-1}.$$
The standard inner product on $\mathbb{R}^{2n}$ is then given by
$$\langle \xi,\eta \rangle =\omega(\xi, J \eta),\qquad \xi, \eta \in \mathbb{R}^{2n}.$$
We suppose that $\Phi$ is a linear symplectomorphisms $\Phi \colon \mathbb{R}^{2n} \to \mathbb{R}^{2n}$, i.e.,
we suppose that $\Phi^* \omega=\omega$. We introduce two Hilbert spaces. Our first Hilbert space is
$$H_1=\Big\{ \xi \in W^{1,2}\big([0,1],\mathbb{R}^{2n}\big): \xi(1)=\Phi \xi(0)\Big\},$$
the space of twisted $W^{1,2}$-loops, and our second Hilbert space is
$$H_0=L^2\big([0,1],\mathbb{R}^{2n}\big).$$
We assume further that for $m \in \mathbb{N}_0$ we have $m$ time-dependent vectors
$$Y_j \in C^\infty\big([0,1],\mathbb{R}^{2n}\big),\qquad 1 \leq j \leq m$$
and an $m\times m$-matrix $A=\{a_{ij}\}_{1\leq i,j \leq m}$, which is symmetric, i.e., $A=A^T$.
We consider the bounded linear operator
$$D \colon H_1 \to H_0$$
which for $\xi \in H_1$ is given by
$$(D\xi)(t)=J\partial_t \xi(t)+\sum_{1 \leq i,j \leq m}a_{ij}\bigg(\int_0^1 \big\langle Y_i(s),\xi(s)\big\rangle ds\bigg) Y_j(t), \qquad
t \in [0,1].$$
\begin{lemma}
The operator $D$ is symmetric with respect to the $L^2$-inner product.
\end{lemma}
\textbf{Proof: } Suppose that $\xi, \eta \in H_1$. We compute
\begin{eqnarray*}
\langle D\xi,\eta\rangle_{L^2}&=&\int_0^1 \big\langle D\xi(t),\eta(t)\big\rangle dt\\
&=&\int_0^1 \big\langle J\partial_t \xi(t),\eta(t)\big \rangle dt\\
& &+\sum_{1 \leq i,j \leq m}a_{ij}\bigg(\int_0^1 \big\langle Y_i(s),\xi(s)\big\rangle ds\bigg) \bigg(\int_0^1 \big \langle Y_j(t),\eta(t)\big\rangle dt\bigg)\\
&=&-\int_0^1 \big\langle J\xi(t),\partial_t\eta(t)\big \rangle dt+\big\langle J\xi(1),\eta(1)\big\rangle-
\big\langle J \xi(0),\eta(0)\big\rangle\\
& &+\sum_{1 \leq i,j \leq m}a_{ij}\bigg(\int_0^1 \big\langle Y_i(s),\eta(s)\big\rangle ds\bigg) \bigg(\int_0^1 \big \langle Y_j(t),\xi(t)\big\rangle dt\bigg)\\
&=&\int_0^1 \big\langle J\partial_t\eta(t),\xi(t)\big \rangle dt+\omega\big(\Phi\xi(0),\Phi\eta(0)\big)-
\omega\big(\xi(0),\eta(0)\big)\\
& &+\sum_{1 \leq i,j \leq m}a_{ij}\bigg(\int_0^1 \big\langle Y_i(s),\eta(s)\big\rangle ds\bigg) \bigg(\int_0^1 \big \langle Y_j(t),\xi(t)\big\rangle dt\bigg)\\
&=&\int_0^1 \big\langle J\partial_t\eta(t),\xi(t)\big \rangle dt+\omega\big(\xi(0),\eta(0)\big)-
\omega\big(\xi(0),\eta(0)\big)\\
& &+\sum_{1 \leq i,j \leq m}a_{ij}\bigg(\int_0^1 \big\langle Y_i(s),\eta(s)\big\rangle ds\bigg) \bigg(\int_0^1 \big \langle Y_j(t),\xi(t)\big\rangle dt\bigg)\\
&=&\int_0^1 \big\langle J\partial_t\eta(t),\xi(t)\big \rangle dt\\
& &+\sum_{1 \leq i,j \leq m}a_{ij}\bigg(\int_0^1 \big\langle Y_i(s),\eta(s)\big\rangle ds\bigg) \bigg(\int_0^1 \big \langle Y_j(t),\xi(t)\big\rangle dt\bigg)\\
&=&\int_0^1 \big\langle D\eta(t),\xi(t)\big\rangle dt\\
&=&\langle D\xi,\eta\rangle_{L^2}.
\end{eqnarray*}
This finishes the proof of the Lemma. \hfill $\square$
\\ \\
The operator $D$ is a compact perturbation of the operator $J\partial_t \colon H_1 \to H_0$. In particular, it
is a Fredholm operator of index zero. Because it is symmetric it becomes an unbounded self-adjoint operator when
interpreted as linear operator $D \colon H_0 \to H_0$ with dense domain $H_1 \subset H_0$.
\begin{lemma}\label{dimest1}
The dimension of the kernel of $D$ can be estimated as follows
$$\mathrm{dim}\big(\ker D\big) \leq 2n+m.$$
\end{lemma}
\textbf{Proof: } Suppose that $\xi \in \mathrm{ker}D$. Then $\xi$ is a solution of the problem
$$J\partial_t \xi(t)+\sum_{1 \leq i,j \leq m}a_{ij}\bigg(\int_0^1 \big\langle Y_i(s),\xi(s)\big\rangle ds\bigg) Y_j(t)=0, \qquad
t \in [0,1],$$
or equivalently
$$\partial_t \xi(t)=\sum_{1 \leq i,j \leq m}a_{ij}\bigg(\int_0^1 \big\langle Y_i(s),\xi(s)\big\rangle ds\bigg) JY_j(t), \qquad
t \in [0,1].$$
Integrating this we obtain
\begin{equation}\label{kerbound}
\xi(t)=\xi(0)+\sum_{1 \leq i,j \leq m}a_{ij}\bigg(\int_0^1 \big\langle Y_i(s),\xi(s)\big\rangle ds\bigg) \int_0^t JY_j(s)ds, \qquad
t \in [0,1].
\end{equation}
Consider the linear map
$$\Gamma \colon \ker D \to \mathbb{R}^{2n} \times \mathbb{R}^m$$
defined for $\xi \in \ker D$ by
$$\Gamma(\xi) = \Bigg(\xi(0),\int_0^1 \big\langle Y_1(s),\xi(s)\big\rangle ds,\ldots, \int_0^1 \big\langle
Y_m(s),\xi(s)\big\rangle ds\Bigg).$$
By (\ref{kerbound}) we see that $\Gamma$ is injective and therefore
$$\mathrm{dim}\big( \ker D\big) \leq \mathrm{dim}\big(\mathbb{R}^{2n} \times \mathbb{R}^m\big)=2n+m.$$
This finishes the proof of the Lemma. \hfill $\square$ 

\section{The commuting case}\label{commcase}

We continue the notation of Section~\ref{selfadjoint}. However, we make some additional assumptions on our data.
Namely we suppose that 
\begin{description}
 \item[(i)] For $1 \leq j \leq m$, the vectors $Y_j$ do not depend on time.
 \item[(ii)] For $1 \leq i,j \leq m$, one has $\omega(Y_i,Y_j)=0$.
\end{description}
Under these stronger hypotheses we obtain in the following Lemma an improvement to Lemma~\ref{dimest1}.
\begin{lemma}\label{dimest2}
Under assumption (i) and (ii) the dimension of the kernel of $D$ can be estimated as follows
$$\mathrm{dim}\big(\ker D\big) \leq 2n.$$
\end{lemma}
\textbf{Proof: } Suppose that $\xi \in \ker D$. Due to assumption (i) the formula (\ref{kerbound}) in the proof of Lemma~\ref{dimest1} 
becomes for $t \in [0,1]$ 
\begin{eqnarray}\label{split1}
\xi(t)&=&\xi(0)+\sum_{1 \leq i,j \leq m}a_{ij}\bigg(\int_0^1 \big\langle Y_i,\xi(s)\big\rangle ds\bigg)tJY_j
\end{eqnarray}
Using (ii) we compute using this formula for $1 \leq i \leq m$
\begin{eqnarray}\label{split2}
\int_0^1 \big\langle Y_i,\xi(s)\big\rangle ds&=&\big\langle Y_i,\xi(0)\big\rangle\\ \nonumber
& &+\sum_{1 \leq i,j \leq m}a_{ij}\bigg(\int_0^1 \big\langle Y_i,\xi(s)\big\rangle ds\bigg)\bigg(\int_0^1 \big\langle Y_i,sJY_j\big\rangle ds\bigg)\\ \nonumber
&=&\big\langle Y_i,\xi(0)\big\rangle\\ \nonumber
& &-\sum_{1 \leq i,j \leq m}a_{ij}\bigg(\int_0^1 \big\langle Y_i,\xi(s)\big\rangle ds\bigg)\bigg(\int_0^1 s\omega(Y_i,Y_j) ds\bigg)\\ \nonumber
&=&\big\langle Y_i,\xi(0)\big\rangle.
\end{eqnarray}
Plugging (\ref{split2}) into (\ref{split1}) we obtain
$$\xi(t)=\xi(0)+t\Bigg(\sum_{1 \leq i,j \leq m}a_{ij}\big\langle Y_i,\xi(0)\big\rangle JY_j\Bigg), \quad t \in [0,1].$$
This shows that $\xi(t)$ is completely determined by $\xi(0) \in \mathbb{R}^{2n}$ and the lemma is proved.  \hfill $\square$

\section{Proof of Theorem\,A}\label{prA}

Since $u \in \mathrm{crit}\mathcal{A}_\mathcal{F}$ we have
$$\partial_t u(t)=X_{df(\overline{H}(u))H}(u(t)), \quad t \in S^1.$$
Abbreviate by $\phi^t_{df(\overline{H}(u))H}$ the flow of the Hamiltonian vector field of ${df(\overline{H}(u))H}$.
With this notion we have
$$u(t)=\phi^t_{df(\overline{H}(u))H} u(0).$$
In particular, since $u$ is periodic, i.e., $u(1)=u(0)$ we can interpret $u(0)$ as a fixed point
$$u(0) \in \mathrm{Fix}\phi^1_{df(\overline{H}(u))H}.$$
Suppose now that
$$\hat{u} \in \ker\big(\mathcal{H}_{\mathcal{A}_\mathcal{F}}(u)\big).$$
For $t \in [0,1]$ abbreviate
$$\xi(t):=d\phi^{-t}_{df(\overline{H}(u))H}\big(u(t)\big)\hat{u}(t) \in T_{u(0)}M.$$
Since $\hat{u}$ is periodic, i.e., $\hat{u}(1)=\hat{u}(0)$, it follows that $\xi \in C^\infty([0,1],T_{u(0)}M)$ satisfies
$$\xi(1)=d\phi^{-1}_{df(\overline{H}(u))H}\big(u(0)\big)\xi(0).$$
Because $\hat{u}$ lies in the kernel of the Hessian of $\mathcal{A}_\mathcal{F}$ at $u$ it follows that $\xi$ is a solution
of the problem
$$\partial_t \xi(t)=d\phi^{-t}_{df(\overline{H}(u))H}X_{d^2f(\overline{H}(u))\big(\int_0^1 dH(u(s))\hat{u}(s)
ds,H\big)}(u(t)),
\quad t\in [0,1].$$
We rewrite this problem as
\begin{eqnarray*}
\partial_t \xi(t)&=&\big(\phi^{t}_{df(\overline{H}(u))H})^*X_{d^2f(\overline{H}(u))\big(\int_0^1 dH(u(s))
d\phi^s_{df(\overline{H}(u))H}\xi(s)ds,H\big)}(u(0))\\
&=&X_{d^2f(\overline{H}(u))\big(\int_0^1 dH \circ \phi^s_{df(\overline{H}(u))H}(u(0))\xi(s)ds,H \circ \phi^t_{df(\overline{H}(u))H}\big)}(u(0)).
\end{eqnarray*}
Choose a basis of $V$ to identify $V=\mathbb{R}^m$ for $m=\dim(V)$. We write $H=(H_1,\ldots, H_m)$ and abbreviate
for $1 \leq i \leq m$ and $t \in \mathbb{R}$ the time-dependent Hamiltonian
$$K_i^t=H_i \circ \phi^t_{df(\overline{H}(u))H}.$$
Setting further
$$a_{ij}:=\frac{\partial^2 f}{\partial x_i \partial x_j}(\overline{H}(u))$$
we can rewrite the above problem as
\begin{eqnarray*}
\partial_t \xi(t)&=&\sum_{1 \leq i,j\leq m} a_{ij} \bigg(\int_0^1 dK^s_i(u(0)) \xi(s) ds\bigg) X_{K_j^t}(u(0))\\
&=&\sum_{1 \leq i,j\leq m} a_{ij} \bigg(\int_0^1 \omega\Big(\xi(s),X_{K_i^s}(u(0))\Big)ds\bigg) X_{K_j^t}(u(0)).
\end{eqnarray*}
Choose further an $\omega_{u(0)}$-compatible complex structure on $T_{u(0)}M$ and abbreviate by
$$\langle \cdot, \cdot \rangle=\omega(\cdot,J\cdot)$$
the corresponding inner product on $T_{u(0)}M$. Using this we can rewrite the above problem equivalently as
$$J \partial_t \xi(t)=-\sum_{1 \leq i,j \leq m}a_{ij}\bigg(\int_0^1 \big \langle J X_{K_i^s}(u(0)),\xi(s)\big \rangle ds \bigg)
J X_{K_j^t}(u(0)).$$
Introducing finally for $t \in [0,1]$ and $1 \leq i \leq m$ the time-dependent vector
$$Y_i(t)=J X_{K_i^t}(u(0)) \in T_{u(0)} M$$
this simplifies to
$$J\partial_t \xi(t)+\sum_{1 \leq i,j \leq m}a_{ij}\bigg(\int_0^1 \big\langle Y_i(s),\xi(s)\big\rangle ds\bigg) Y_j(t)=0, \qquad t \in [0,1].$$
Now Theorem\,A follows from Lemma~\ref{dimest1}.

\section{Proof of Theorem\,B}\label{prB}

We continue the notion of Section~\ref{prA}. Since the pair $\mathcal{F}=(f,H)$ is commuting we obtain that
$$K_i^t=H_i \circ \phi^t_{df(\overline{H}(u))H}=H_i$$
for $1 \leq i \leq m$ so that
$$Y_i(t)=J X_{K_i^t}(u(0)) \in T_{u(0)} M=JX_{H_i}(u(0))$$
does not depend on time. Using again that $\mathcal{F}$ is commuting we get that
$$\omega(X_{H_i},X_{H_j})=0$$
for $1 \leq i,j \leq m$, so that since $J$ is $\omega$-compatible
$$\omega(Y_i,Y_j)=\omega(JX_{H_i},JX_{H_j})(u(0))=\omega(X_{H_i},X_{H_j})(u(0))=0.$$
This shows that the vectors $Y_i$ meet the requirements (i) and (ii) in Section~\ref{commcase} and therefore Theorem\,B
follows from Lemma~\ref{dimest2}.

\appendix

\section{Symmetries}\label{symmetries}

In this appendix we explain how critical points of the action functional $\mathcal{A}_\mathcal{F}$ can be iterated. 
\\ \\
We denote by $\mathbb{Z}^*$ the monoid which as a set is given by
$$\mathbb{Z}^*=\mathbb{Z} \setminus \{0\}$$
and product given by the usual multiplication in $\mathbb{Z}$. We further consider the monoid $\mathbb{Z}^* \ltimes S^1$
given by the semidirect product of the monoid $\mathbb{Z}^*$ and the group $S^1$ whose multiplication is defined by
$$(n_1,r_1)(n_2,r_2)=(n_1 n_2,n_1 r_2+r_1), \qquad (n_1,r_1),(n_2,r_2) \in \mathbb{Z}^* \ltimes S^1.$$
The monoid $\mathbb{Z}^* \ltimes S^1$ acts on the component of contractible loops $\mathcal{L}$ of the free loop space 
for $u \in \mathcal{L}$ and $(n,r) \in \mathbb{Z}^* \ltimes S^1$ by
$$(n,r)_* u(t)=u\big(n(t+r)\big), \qquad t \in S^1.$$
We consider now a pair $\mathcal{F}=(f,H)$, where $H \colon M \to V$ is a smooth map from $M$ to a finite dimensional
real vector space $V$ and $f \colon W \to \mathbb{R}$ is a smooth map defined on an open subset $W$ of $V$. Note
that $\overline{H} \colon \mathcal{L} \to V,\,\,u \mapsto \int_0^1 H(u(t))dt$
is invariant under the action of the monoid $\mathbb{Z}^* \ltimes S^1$ so that $\mathcal{L}_\mathcal{F}=\overline{H}^{-1}(W)$ is invariant as well. For $(n,r) \in \mathbb{Z} \ltimes S^1$ we define the pullback of the tuple $\mathcal{F}$ by
$$(n,r)^*(f,H)=\big(\tfrac{1}{n}f,H\big).$$
Observe that we have
$$\mathcal{L}_\mathcal{F}=\mathcal{L}_{(n,r)^*\mathcal{F}}.$$
The pullback of the action functional $\mathcal{A}_\mathcal{F} \colon \mathcal{L}_\mathcal{F} \to \mathbb{R}$ is given
by
$$(n,r)^*\mathcal{A}_\mathcal{F}=\mathcal{A}_{(n,r)^*\mathcal{F}}.$$
In particular, if $u \in \mathcal{L}_\mathcal{F}$ and $(n,r)_* u \in \mathrm{crit}\mathcal{A}_\mathcal{F}$ it follows that
$u$ is a critical point of $\mathcal{A}_{(n,r)^*\mathcal{F}}$. In view of the critical point equation in Lemma~\ref{crit}
the converse is true as well so that we have
\begin{prop}
Suppose that $u \in \mathcal{L}_\mathcal{F}$ and $(n,r) \in \mathbb{Z}^* \ltimes S^1$. Then $u$ is a critical point of
$\mathcal{A}_{(n,r)^*\mathcal{F}}$ if and only if $(n,r)_* u$ is a critical point of $\mathcal{A}_\mathcal{F}$.
\end{prop}

\end{document}